\documentclass{amsart}
\usepackage{amsmath}
\usepackage{amsthm}
\usepackage{amssymb}
\usepackage{epsfig}

\theoremstyle{plain}
\newtheorem{thm}{Theorem}[section]
\newtheorem*{unnum thm}{Theorem}

\newtheorem{prop}[thm]{Proposition}
\newtheorem{cor}[thm]{Corollary}

\theoremstyle{definition}
\newtheorem{example}[thm]{Example}
\newtheorem{defn}[thm]{Definition}
\theoremstyle{remark}

\newcommand{\eqdef}{\ensuremath :=}
\newcommand{\bm}[1]{\mbox{\boldmath $#1$}}
\newcommand{\ul}[1]{\underline{#1}}
\newcommand{\B}{\mathfrak{S}_n^B}
\newcommand{\D}{\mathfrak{S}_n^D}

\title{Pattern Avoidance and the Bruhat Order}
\author{Bridget Eileen Tenner}
\address{Department of Mathematics, Massachusetts Institute of Technology, 77 Massachusetts Ave, Cambridge, MA 02139, USA}
\email{bridget@math.mit.edu}
\date{April 13, 2006}

\begin{document}

\begin{abstract}

The structure of order ideals in the Bruhat order for the symmetric group is elucidated via permutation patterns.  A method for determining non-isomorphic principal order ideals is described and applied for small lengths.  The permutations with boolean principal order ideals are characterized.  These form an order ideal which is a simplicial poset, and its rank generating function is computed.  Moreover, the permutations whose principal order ideals have a form related to boolean posets are also completely described.  It is determined when the set of permutations avoiding 
a particular set of patterns is an order ideal, and the rank generating functions of these ideals are computed.  Finally, the Bruhat order in types $B$ and $D$ is studied, and the elements with boolean principal order ideals are characterized and enumerated by length.

\end{abstract}

\maketitle

\section{Introduction}\label{section intro}

This paper studies the interplay between the Bruhat order and permutation patterns, with particular emphasis on these relationships in the symmetric group.  The principal order ideals in particular are considered, and several results are described which emphasize the relationship between permutation patterns and reduced decompositions.  The final section of the paper discusses the Bruhat order for types $B$ and $D$, although in less depth than the type $A$ discussion.

The finite Coxeter groups of types $A$, $B$, and $D$ have combinatorial interpretations as permutations, signed permutations, and signed permutations with certain restrictions.  The combinatorial aspects of Coxeter groups are treated in \cite{bjornerbrenti}.  Although these groups are classical objects with a bountiful literature, there are still many open questions, particularly in reference to patterns and the Bruhat order.

Following the work of Simion and Schmidt in \cite{simion}, there has been a surge of interest in permutation patterns.  Although many intriguing results have been shown, some of the most basic questions remain unanswered.  However, recent work (see \cite{rdpp}) has uncovered connections between reduced decompositions and permutation patterns that may prove useful to resolving some of these issues.

The Bruhat order is a partial ordering of Coxeter group elements, and it plays a remarkably significant role in the study of these groups.  Somewhat surprisingly, very little is known about its structure, particularly in terms of its order ideals and intervals.  The results presented here elucidate some pattern-related facts about this structure.  When combined with the relationship between reduced decompositions and patterns in \cite{rdpp}, these are significant steps towards understanding the more general structural aspects of this partial order.  This paper primarily examines the structure of the Bruhat order of the symmetric group.

Subsequent to Sections~\ref{defs section} and~\ref{bruhat order section} which introduce the concepts that will be discussed throughout the paper, Section~\ref{isomorphism section} considers isomorphism classes of order ideals in the Bruhat order of the symmetric group.  The notion of a \emph{decomposable} order ideal is introduced, and this provides a method for describing the non-isomorphic principals order ideals of a given length.

Section~\ref{boolean section} classifies all permutations with boolean principal order ideals.  As shown in Theorem~\ref{boolean thm}, these are exactly those permutations that avoid the patterns $321$ and $3412$.  The permutations with this property are enumerated by length in Corollary~\ref{boolean by length cor}.  Additionally, permutations with ``nearly boolean'' principal order ideals are discussed, along with the size and description of their ideals.

A more general classification occurs in Section~\ref{power section}.  There the permutations with principal order ideals isomorphic to a power of $B(w_0^{(k)})$, for $k \ge 3$, are entirely classified as those in which every inversion is in exactly one decreasing subsequence of length $k$.  Although this characterization (Theorem~\ref{power thm}) is again stated in terms of patterns, it is not exactly pattern avoidance.

Section~\ref{pattern ideal section} examines sets of permutations avoiding either one or two patterns, and determines exactly when these sets are order ideals in the Bruhat order.  This property holds in only a few situations, each of which can be enumerated by length.

Expanding on the results of Section~\ref{boolean section}, Section~\ref{bd bruhat section} examines the Bruhat order for the finite Coxeter groups of types $B$ and $D$.  In particular, those elements with boolean principal order ideals are defined in Theorems~\ref{type b boolean} and~\ref{type d boolean}.  Once again, permutation patterns emerge, although now for signed permutations, and the avoidance of certain patterns is equivalent to having a boolean principal order ideal.  While the case for type $A$ required avoiding only two patterns, it is necessary to avoid ten patterns in type $B$, and twenty patterns must be avoided to have a boolean principal order ideal in type $D$.  For types $B$ and $D$, the elements avoiding these patterns are enumerated by length in Corollaries~\ref{boolean b by length} and~\ref{boolean d by length}.

\section{Definitions and background}\label{defs section}

Let $\mathfrak{S}_n$ be the group of permutations on $n$ elements, and let $[n]$ denote the set of integers $\{1, \ldots, n\}$.  An element $w \in \mathfrak{S}_n$ is the bijection on $[n]$ mapping $i \mapsto w(i)$.  A permutation will be written in one-line notation as $w = w(1)w(2) \cdots w(n)$.

\begin{example}
$4213 \in \mathfrak{S}_4$ maps $1 \mapsto 4$, $2 \mapsto 2$, $3 \mapsto 1$, and $4 \mapsto 3$.
\end{example}

\begin{defn}
An \emph{inversion} in $w$ is a pair $(i,j)$ such that $i<j$ and $w(i)>w(j)$.
\end{defn}

Let $[\pm n]$ denote the set of integers $\{\pm 1, \ldots, \pm n\}$.  For ease of notation, a negative sign may be written beneath an integer: $\ul{i} \eqdef -i$.

\begin{defn}
A \emph{signed permutation} of $[\pm n]$ is a bijection $w$ with the requirement that $w(\ul{i}) = \ul{w(i)}$.  Let $\B$ denote the signed permutations of $[\pm n]$.
\end{defn}

An element $w \in \B$ is entirely defined by $w(1), \ldots, w(n)$.  Therefore one-line notation will again be used, although some values may now be negative.

\begin{example}\label{type B example}
$\ul{4}21\ul{3} \in \mathfrak{S}_4^B$ maps $\pm 1 \mapsto \mp 4$, $\pm 2 \mapsto \pm 2$, $\pm 3 \mapsto \pm 1$, and $\pm 4 \mapsto \mp 3$.
\end{example}

The Coxeter groups studied here are the finite Coxeter groups of types $A$, $B$, and $D$.  More thorough discussions of general Coxeter groups appear in \cite{bjornerbrenti} and \cite{humphreys}.

Define involutions, called \emph{simple reflections}, on $[\pm n]$ as follows:
\begin{eqnarray*}
s_i&=&1\cdots (i-1) (i+1)i (i+2) \cdots n \text{ for } i \in [n-1];\\
s_0&=&\ul{1}2\cdots n; \text{ and}\\
s_{1'}&\eqdef&s_0s_1s_0 = \ul{21}3\cdots n.
\end{eqnarray*}
\noindent These definitions indicate that the following braid relations hold:
\begin{eqnarray}
\label{relation ij commute}s_is_j &=& s_js_i \text{ for } i,j \in [0,n-1]\text{ and }|i-j|>1;\\
\label{relation hati commute}s_{1'}s_i &=& s_is_{1'} \text{ for } i \in \big([n-1] \setminus \{2\}\big);\\
\label{relation ii+1}s_is_{i+1}s_i &=& s_{i+1}s_is_{i+1} \text{ for } i \in [n-2];\\
\label{relation 01}s_0s_1s_0s_1 &=& s_1s_0s_1s_0; \text{ and}\\
\label{relation hat2}s_{1'}s_2s_{1'} &=& s_2s_{1'}s_2.
\end{eqnarray}

\begin{defn}\label{type A def}
The finite Coxeter group of \emph{type $A$} is the symmetric group $\mathfrak{S}_n$, for some $n$.  This group is generated by $\{s_1, \ldots, s_{n-1}\}$.
\end{defn}

\begin{defn}\label{type B def}
The finite Coxeter group of \emph{type $B$} is the hyperoctahedral group $\B$, for some $n$.  This group is generated by $\{s_0, s_1, \ldots, s_{n-1}\}$.
\end{defn}

\begin{defn}\label{type D def}
The finite Coxeter group of \emph{type $D$} is the subgroup $\D$ of $\B$ consisting of signed permutations whose one-line notation contains an even number of negative values.  This group is generated by $\{s_{1'}, s_1, \ldots, s_{n-1}\}$.
\end{defn}

\begin{defn}\label{reduced decomposition def}
Let $W$ be a Coxeter group generated by the simple reflections $\mathcal{S}$.  For $w \in W$, if $w = s_{i_1} \cdots s_{i_{\ell}}$ and $\ell$ is minimal among all such expressions, then the string $i_1 \cdots i_{\ell}$ is a \emph{reduced decomposition} of $w$ and $\ell$ is the \emph{length} of $w$, denoted $\ell(w)$.  The set $R(w)$ consists of all reduced decompositions of $w$.
\end{defn}

The permutation $w_0 \eqdef n\cdots 21 \in \mathfrak{S}_n$, is aptly named the \emph{longest element} in $\mathfrak{S}_n$, and $\ell(w_0) = \binom{n}{2}$.  If $n$ is unclear from the context, this may be denoted $w_0^{(n)}$.  Analogously, the \emph{longest element} in $\B$ is $\ul{12}\cdots \ul{n}$, and the longest element in $\D$ is $\ul{12} \cdots \ul{n}$ if $n$ is even, and $1\ul{23}\cdots\ul{n}$ if $n$ is odd.

\begin{defn}
A consecutive substring of a reduced decomposition is a \emph{factor}.
\end{defn}

Multiplication here follows the standard that a function is to the left of its input.  Thus, if $i \in [n-1]$, then $s_iw$ interchanges the positions of the values $i$ and $i+1$ (and $\ul{i}$ and $\ul{i+1}$) in $w$, whereas $ws_i = w(1) \cdots w(i+1)w(i) \cdots w(n)$.

The classical notion of (unsigned) permutation pattern avoidance is as follows.

\begin{defn}\label{defn:patterns}
Fix $w \in \mathfrak{S}_n$ and $p \in \mathfrak{S}_k$ for $k \le n$.  The permutation $w$ \emph{contains the pattern $p$}, or \emph{contains a $p$-pattern}, if there exist $i_1 < \cdots < i_k$ such that $w(i_1) \cdots w(i_k)$ is in the same relative order as $p(1) \cdots p(k)$.  That is, $w(i_h) < w(i_j)$ if and only if $p(h) < p(j)$.  If $w$ does not contain $p$, then $w$ \emph{avoids} $p$, or is \emph{$p$-avoiding}.
\end{defn}

In \cite{rdpp}, reduced decompositions are analyzed in conjunction with pattern containment in $\mathfrak{S}_n$.  This coordinated approach yields a number of significant results, and the main theorem of \cite{rdpp} is the vexillary characterization below.

\begin{defn}
The \emph{shift} of a string $\bm{i} = i_1 \cdots i_\ell$ by $M \in \mathbb{N}$ is
\begin{equation*}
\bm{i}^M := (i_1 + M) \cdots (i_\ell + M).
\end{equation*}
\end{defn}

\begin{unnum thm}[Vexillary characterization]
The permutation $p$ is vexillary if and only if, for every permutation $w$ containing a $p$-pattern, there exists a reduced decomposition $\bm{j} \in R(w)$ containing a shift of some $\bm{i} \in R(p)$ as a factor.
\end{unnum thm}

The definition of patterns in signed permutations requires an extra clause.

\begin{defn}\label{defn:signed patterns}
Fix $w \in \B$ and $p \in \mathfrak{S}_k^B$ for $k \le n$.  The permutation $w$ \emph{contains the pattern $p$} if there exist $0 < i_1 < \cdots < i_k$ such that
\begin{enumerate}
\item $w(i_j)$ and $p(j)$ have the same sign; and 
\item $|w(i_1)| \cdots |w(i_k)|$ is in the same relative order as $|p(1)| \cdots |p(k)|$.
\end{enumerate}
If $w$ does not contain $p$, then $w$ \emph{avoids} $p$, or is \emph{$p$-avoiding}.
\end{defn}

\begin{example}
Let $w = \ul{4}21\ul{3}$, $p = \ul{3}1\ul{2}$, $q = 31\ul{2}$, and $r = 132$.  Then $\ul{4}1\ul{3}$ and $\ul{4}2\ul{3}$ are both occurrences of $p$ in $w$.  The signed permutation $w$ is $q$- and $r$-avoiding.
\end{example}

If $w$ has a $p$-pattern, with $\{i_1, \ldots, i_k\}$ as in Definitions~\ref{defn:patterns} and~\ref{defn:signed patterns}, then $w(i_1) \cdots w(i_k)$ is an \emph{occurrence} of $p$ in $w$.  Define $\langle p(j) p(j+1) \cdots p(j+m) \rangle$ to be $w(i_j) w(i_{j+1}) \cdots w(i_{j+m})$.  Occurrences of $p$ will be distinguished by subscripts: $\langle \ \rangle_i$.

\begin{example}
Let $w = 7413625$, $p = 1243$, and $q = 1234$.  Then $1365$ is an occurrence of $p$, with $\langle 1 \rangle = 1$, $\langle 2 4 \rangle = 36$, and $\langle 3 \rangle = 5$.  The permutation $w$ avoids $q$.
\end{example}

\section{Bruhat order}\label{bruhat order section}

Standard terminology from the theory of partially ordered sets will be used throughout this paper.  Good sources for information on this topic are \cite{ec1} and \cite{trotter}.

The Bruhat order is a partial ordering that can be placed on a Coxeter group.

\begin{defn}\label{bruhat order def}
Fix a Coxeter group $W$ generated by the simple reflections $\mathcal{S}$.  Let $\overline{\mathcal{S}} = \{wsw^{-1}: w \in W, s \in \mathcal{S}\}$.  For $w, w' \in W$, write $w \lessdot w'$ if $\ell(w') = \ell(w) + 1$ and $w' = \overline{s}w$ for some $\overline{s} \in \overline{\mathcal{S}}$.  The covering relations $\lessdot$ give the \emph{(strong) Bruhat order}.
\end{defn}

As discussed in \cite{bjornerbrenti}, the Bruhat order does not favor left or right multiplication, despite the appearance of Definition~\ref{bruhat order def}.

\begin{example}
$7314625 \gtrdot 7312645$.  Although not a covering relation, $7314625 > 1374625$ in the Bruhat order.
\end{example}

This partial order has many properties which are discussed and proved in \cite{bjornerbrenti}.  The property most relevant to this work, the subword property, gives an equivalent definition of the Bruhat order in terms of reduced decompositions.

\begin{unnum thm}[Subword property]
Let $W$ be a Coxeter group and $w,w' \in W$.  Choose $i_1 \cdots i_{\ell'} \in R(w')$.  Then $w \le w'$ in the Bruhat order if and only if there exists $j_1 \cdots j_{\ell} \in R(w)$ which is a subword of $i_1 \cdots i_{\ell'}$.
\end{unnum thm}

\begin{figure}[htbp]
\centering
\epsfig{file=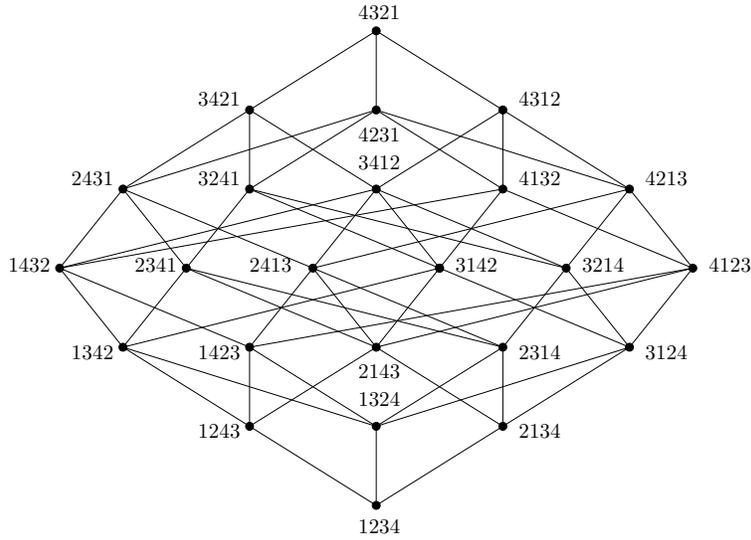, width=4in}
\caption{The Bruhat order for $\mathfrak{S}_4$.}\label{bruhat4321 fig}
\end{figure}

\begin{figure}[htbp]
\centering
\epsfig{file=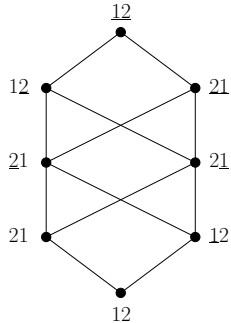, width=1.16in}
\caption{The Bruhat order for $\mathfrak{S}_2^B$.}\label{bruhat-1-2.fig}
\end{figure}

The Bruhat order for a Coxeter group gives a graded poset where the rank function is the length of an element.  Figures~\ref{bruhat4321 fig} and~\ref{bruhat-1-2.fig} give the Hasse diagrams for the Bruhat order on $\mathfrak{S}_4$ and $\mathfrak{S}_2^B$, respectively.

Other properties of the Bruhat order include that it is an Eulerian poset (shown by Verma in \cite{verma}) and that it is CL-shellable (shown by Bj\"{o}rner and Wachs in \cite{bjornerwachs}).  Bj\"{o}rner and Wachs also show that every open interval in the poset is topologically a sphere.  For more information, see \cite{bjornerwachs}, \cite{ec1}, and \cite{verma}.

Understanding intervals in the Bruhat order is made substantially simpler by Dyer's result in \cite{dyer}: for any $\ell$, there are only finitely many non-isomorphic intervals of length $\ell$ in the Bruhat order of finite Coxeter groups.  The length $4$ intervals have been classified, as have length $5$ intervals in $\mathfrak{S}_n$, by Hultman in~\cite{hultman1} and~\cite{hultman2}.

Order ideals, specifically principal order ideals, are studied here.

\begin{defn}
Let $W$ be a finite Coxeter group.  For $w \in W$, let
\begin{equation*}
B(w) = \{v \in W : v \le w\}
\end{equation*}
\noindent be the principal order ideal of $w$ in the Bruhat order for $W$.
\end{defn}

In \cite{sjostrand}, Sj\"{o}strand studies $B(w)$ in relation to rook configurations and Ferrers boards.  He also gives a polynomial time recurrence for computing some $|B(w)|$.

\section{Isomorphism classes of principal order ideals}\label{isomorphism section}

This section examines isomorphism classes of principal order ideals in the Bruhat order for the symmetric group.

\begin{defn}
A permutation $w \in \mathfrak{S}_n$ is \emph{decomposable} if $B(w) \cong B(u) \times B(v)$ for some $u \in \mathfrak{S}_m$ and $v \in \mathfrak{S}_{m'}$ where $m,m' < n$.  Otherwise it is \emph{indecomposable}.
\end{defn}

The following results are straightforward to show, and the proofs are omitted.

\begin{prop}\label{concat prop}
A permutation $w \in \mathfrak{S}_n$ is decomposable if and only if there exist $M \in [n-2]$ and $\bm{g}\bm{h}$ or $\bm{h}\bm{g}$ in $R(w)$, such that $\bm{g}$ is comprised of letters in $[M]$ and $\bm{h}$ is comprised of letters in $[M+1,n-1]$.  Equivalently, $w$ is indecomposable if and only if there is a substring $M(M+1)M$ or $(M+1)M(M+1)$ in every element of $R(w)$, for all $M \in [n-2]$.
\end{prop}

Suppose that $w$ is decomposable, and keep the notation of Proposition~\ref{concat prop}.  Let $\bm{H}$ be such that $\bm{H}^M = \bm{h}$, and let $u \in \mathfrak{S}_{M+1}$ and $v \in \mathfrak{S}_{n-M}$ be such that $\bm{g} \in R(u)$ and $\bm{H} \in R(v)$.  Then $B(w) \cong B(u) \times B(v)$. 

\begin{cor}
The decomposable permutations form an order ideal in $\mathfrak{S}_n$.
\end{cor}

These results greatly simplify the problem of determining the non-isomorphic principal order ideals in the Bruhat order of the symmetric group, as indicated in Table~\ref{non-iso table}.  Entries for greater lengths can be similarly deduced.

\begin{table}[htbp]
\centering
\begin{tabular}{|l|cccccc|}
\hline
\rule[-2mm]{0mm}{6mm}Length& $ 0$ & $1$ & $2$ & $3$ & $4$ & $5$\\
\hline
\rule[0mm]{0mm}{4mm}Reduced& $\emptyset$ & $1$ & $12$ & $123$ & $1234$ & $12345$\\
decompositions& & & & $121$ & $1214$ & $12146$\\
& & & & & $2132$ & $21325$\\
& & & & & & $12321$\\
& & & & & & $21232$\\
\hline
\end{tabular}
\caption{Representative reduced decompositions giving all non-isomorphic principal order ideals of length at most $5$ in $\mathfrak{S}_n$.}\label{non-iso table}
\end{table}

\begin{figure}[htbp]
\centering
\epsfig{file=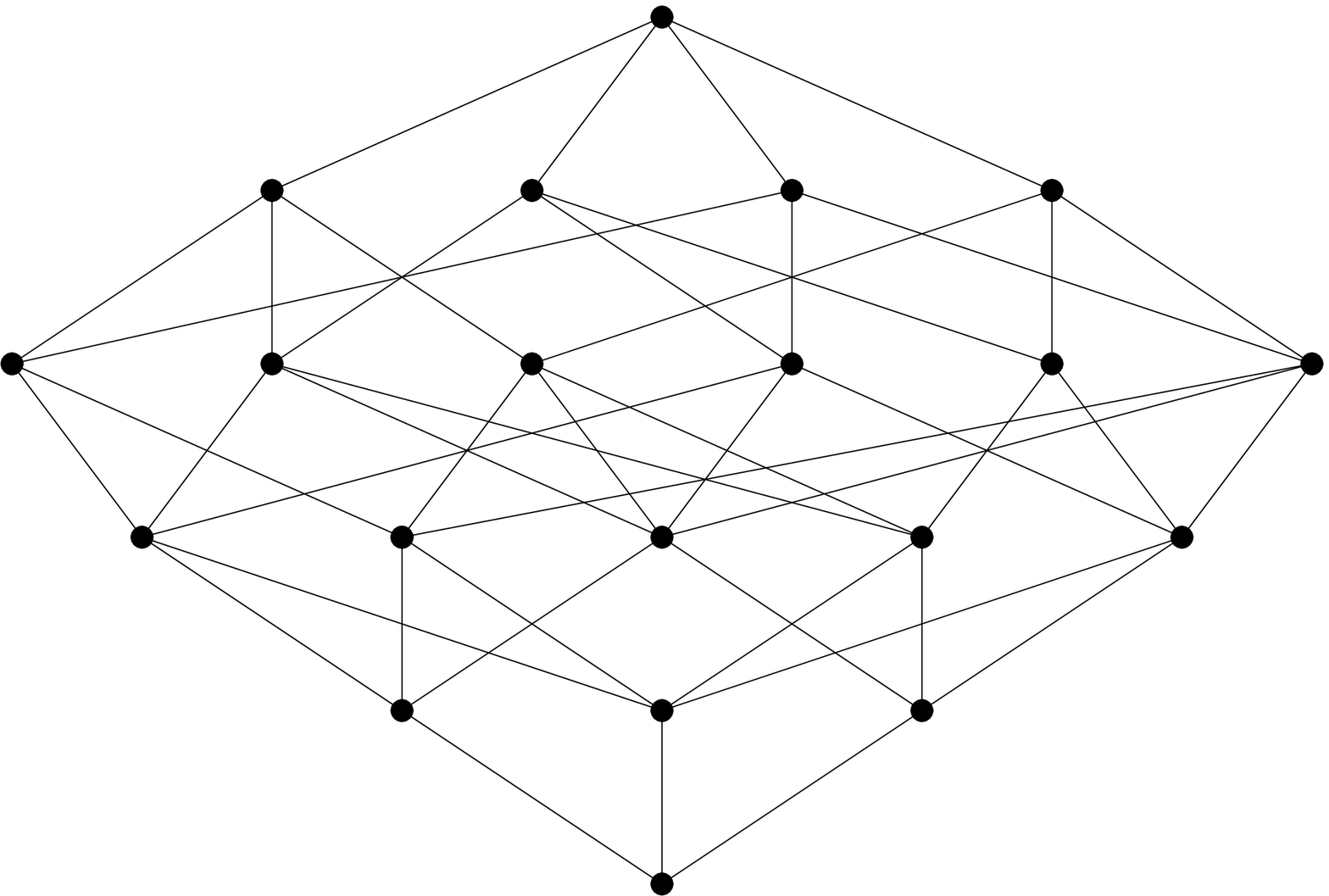, width=2.62in}
\caption{The principal order ideal $B(4231)$.}\label{B(4231) fig}
\end{figure}

\begin{figure}[htbp]
\centering
\epsfig{file=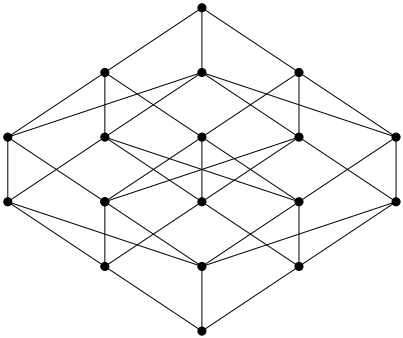, width=2.08in}
\caption{The principal order ideal $B(3421)$.}\label{B(3421) fig}
\end{figure}

The entries in Table~\ref{non-iso table} correspond to principal order ideals that can be described by only a few posets, many of which are depicted in Figures~\ref{B(4231) fig}, \ref{B(3421) fig}, \ref{B(321) fig}, and~\ref{B(3412) fig}.  Table~\ref{non-iso ideals} indicates these principal order ideals.

\begin{table}[htbp]
\centering
\begin{tabular}{|l|cccccc|}
\hline
\rule[-2mm]{0mm}{6mm}Length& $ 0$ & $1$ & $2$ & $3$ & $4$ & $5$\\
\hline
\rule[0mm]{0mm}{4mm}Principal& $B(1)$ & $B(21)$ & $B(21)^2$ & $B(21)^3$ & $B(21)^4$ & $B(21)^5$\\
order& & & & $B(321)$ & $B(321) \times B(21)$ & $B(321) \times B(21)^2$\\
ideals& & & & & $B(3412)$ & $B(3412) \times B(21)$\\
& & & & & & $B(4231)$\\
& & & & & & $B(3421)$\\
\hline
\end{tabular}
\caption[The non-isomorphic principal order ideals of length at most $5$ in $\mathfrak{S}_n$.]{The non-isomorphic principal order ideals of length at most $5$ in $\mathfrak{S}_n$.  The entries correspond to those in Table~\ref{non-iso table}.}\label{non-iso ideals}
\end{table}

Not all intervals in $\mathfrak{S}_n$ can be principal order ideals.  The discrepancies in the number of non-isomorphic intervals that can appear and the number of non-isomorphic principal order ideals that can appear are displayed in Table~\ref{discrepancies table}.  The theorem that there are only finitely many non-isomorphic intervals of a given length in the Bruhat order is due to Dyer in \cite{dyer}, and the quantitative results for small lengths are due to Jantzen (see \cite{jantzen}) and Hultman (see \cite{hultman1} and \cite{hultman2}).

\begin{table}[htbp]
\centering
\begin{tabular}{|r|c|c|c|c|c|c|}
\hline
\rule[-2mm]{0mm}{6mm}Length & $ 0$ & $1$ & $2$ & $3$ & $4$ & $5$\\
\hline
\rule[-2mm]{0mm}{6mm}$\#$ Non-isomorphic intervals & 1 & 1 & 1 & 3 & 7 & 25\\
\hline
\rule[-2mm]{0mm}{6mm}$\#$ Non-isomorphic $B(w)$ & 1 & 1 & 1 & 2 & 3 & 5\\
\hline
\end{tabular}
\caption{Comparing the non-isomorphic intervals and the non-isomorphic principal order ideals of length at most $5$ in $\mathfrak{S}_n$.}\label{discrepancies table}
\end{table}

\section{Boolean principal order ideals}\label{boolean section}

\begin{defn}\label{boolean def}
The \emph{boolean poset} $B_r$ is the set of subsets of $[r]$ ordered by set inclusion.  A poset is \emph{boolean} if it is isomorphic to $B_r$ for some $r$.
\end{defn}

\begin{figure}[htbp]
\centering
\epsfig{file=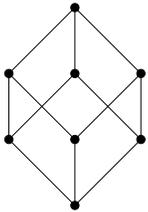, width=.74in}
\caption{The boolean poset $B_3$.}\label{boolean B_3 fig}
\end{figure}

Because $B_1 \cong B(21)$, a poset is boolean if and only if it is isomorphic to $B(21)^r$ for some $r$.  For example, the poset depicted in Figure~\ref{boolean B_3 fig} is isomorphic to $B(21)^3$.  The goal of this section is to determine exactly when the principal order ideal $B(w)$ is boolean, for $w \in \mathfrak{S}_n$.

\begin{defn}
Let $w$ be a permutation in $\mathfrak{S}_n$.  If the poset $B(w)$ is boolean, then $w$ is a \emph{boolean} permutation.
\end{defn}

\begin{thm}\label{boolean thm}
The permutation $w$ is boolean if and only if $w$ is $321$- and $3412$-avoiding.
\end{thm}

\begin{proof}
The poset $B(w)$ is graded of rank $\ell \eqdef \ell(w)$.  Thus, if $B(w)$ is boolean, then $B(w) \cong B_{\ell}$.  Fix $\bm{i} = i_1 \cdots i_{\ell} \in R(w)$.  By the subword property, it must be possible to delete any subset of $\{i_1, \ldots, i_{\ell}\}$ and obtain a string which is still reduced.  From this, it is straightforward to show that boolean permutations are exactly those which have a reduced decomposition with no repeated letters.  In fact, if one reduced decomposition has this property, then all reduced decompositions do.

Suppose that $w$ is not boolean.  Then there is a reduced decomposition of $w$ with a repeated letter.  Therefore, there exists a reduced decomposition of $w$ with one of the following factors, for some $M$.
\begin{eqnarray}
\label{121 factor}&(121)^M&\\
\label{2132 factor}&(2132)^M&
\end{eqnarray}

\noindent By the vexillary characterization of \cite{rdpp}, a factor as in equation~\eqref{121 factor} indicates that $w$ has a $321$-pattern.  Similarly, a factor as in equation~\eqref{2132 factor} implies either a $3412$-pattern, or a $4312$-, $3421$-, or $4321$-pattern.  The latter three all contain the pattern $321$, so a repeated letter in a reduced decomposition implies that $w$ has a $321$- or a $3412$-pattern.

Conversely, $321$ and $3412$ are vexillary.  Thus, by the vexillary characterization, if either pattern appears then a reduced decomposition has a repeated letter.
\end{proof}

Boolean permutations were previously enumerated by West in \cite{west} and Fan in \cite{fan}, although under different guises.

\begin{cor}[Fan, West]\label{number boolean}
The number of boolean permutations in $\mathfrak{S}_n$ is $F_{2n-1}$, where $\{F_0, F_1, \ldots\}$ are the Fibonacci numbers.  This is sequence A001519 in \cite{oeis}.
\end{cor}

The boolean permutations can also be enumerated by length.

\begin{cor}\label{boolean by length cor}
Let $L(n,k) \eqdef \# \{w \in \mathfrak{S}_n : \ell(w) = k \text{ and } w \text{ is boolean}\}$.  Then
\begin{equation}\label{L(n,k) equation}
L(n,k) = \sum_{i=1}^k \binom{n-i}{k+1-i} \binom{k-1}{i-1},
\end{equation}
\noindent where the (empty) sum for $k=0$ is defined to be $1$.
\end{cor}

\begin{proof}
The result is proved by induction.  First, observe that there is exactly one permutation in $\mathfrak{S}_n$ of length $0$, and it is boolean.  There are $n-1$ permutations in $\mathfrak{S}_n$ of length $1$, and these are all boolean.  Letting $k=1$ in equation~\eqref{L(n,k) equation} yields
\begin{equation*}
\sum_{i=1}^1 \binom{n-i}{1+1-i} \binom{1-1}{i-1} = \binom{n-1}{1}\binom{0}{0} = n-1,
\end{equation*}
\noindent so the corollary holds for $k \le 1$ and any $n>k$.

Assume the result for all $k \in [0,K)$ and $n \in [1,N)$.  A boolean permutation avoids the patterns $321$ and $3412$.  Suppose $w \in \mathfrak{S}_N$ is a boolean permutation with $\ell(w) = K$, and consider the location of $N$ in the one-line notation of $w$.
\begin{itemize}
\item If $w(N) = N$, then $w(1) \cdots w(N-1) \in \mathfrak{S}_{N-1}$ can be any boolean permutation of length $K$.
\item If $w(N-1) = N$, then $w(1) \cdots w(N-2) w(N) \in \mathfrak{S}_{N-1}$ can be any boolean permutation of length $K-1$.
\item If $w(N-2) = N$, then $w(N) = N-1$.  Thus, $w(1) \cdots w(N-3) w(N-1) \in \mathfrak{S}_{N-2}$ can be any boolean permutation of length $K-2$.
\item If $w(N-3) = N$, then $w(N) = N-1$ and $w(N-1) = N-2$.  Therefore $w(1) \cdots w(N-4) w(N-2) \in \mathfrak{S}_{N-3}$ can be any boolean permutation of length $K-3$.
\item $\ldots$
\end{itemize}
\noindent Thus $L(N,K) = L(N-1,K) + \sum\limits_{i=1}^K L(N-i,K-i)$, and
\begin{equation*}
L(N,K) = L(N-1,K-1) + L(N-1,K-1) + L(N-2,K-1) + \cdots + L(K,K-1).
\end{equation*}

By the inductive assumptions and basic facts about binomial coefficients,
\begin{eqnarray*}
L(N,K) &=& \sum_{i=1}^{K-1} \binom{N-1-i}{K-i}\binom{K-2}{i-1} + \sum_{j=K}^{N-1} \sum_{i=1}^{K-1} \binom{j-i}{K-i}\binom{K-2}{i-1}\\
&=& \sum_{i=1}^{K-1} \binom{N-1-i}{K-i}\binom{K-2}{i-1} + \sum_{i=1}^{K-1} \binom{N-i}{K+1-i}\binom{K-2}{i-1}\\
&=& \sum_{i=2}^{K}\binom{N-i}{K+1-i}\binom{K-2}{i-2} + \sum_{i=1}^{K-1} \binom{N-i}{K+1-i}\binom{K-2}{i-1}\\
&=& \binom{N-K}{1}\binom{K-2}{K-2} + \sum_{i=2}^{K-1} \binom{N-i}{K+1-i} \left(\binom{K-2}{i-2} + \binom{K-2}{i-1}\right)\\
& & \hspace{1in} + \binom{N-1}{K}\binom{K-2}{0}\\
&=& \sum_{i=1}^K \binom{N-i}{K+1-i} \binom{K-1}{i-1}.
\end{eqnarray*}
\end{proof}

The numbers $L(n,k)$ are equal to the numbers $T(n,n-k)$ in sequence A105306 of \cite{oeis}.  From this, it is straightforward to compute the generating function
\begin{equation}\label{boolean gen fxn}
\sum_{n,k}L(n,k)t^kz^n = \frac{z(1-zt)}{1-2zt-z(1-zt)}.
\end{equation}

For small $n$ and $k$, the values $L(n,k)$ are displayed in Table~\ref{L(n,k) table}.  

\begin{table}[htbp]
\centering
\begin{tabular}{|c|cccccccc|}
\hline
\rule[-2mm]{0mm}{6mm}$L(n,k)$ & $k=0$ & $1$ & $2$ & $3$ & $4$ & $5$ & $6$ & $7$\\
\hline
\rule[0mm]{0mm}{4mm}$n=1$ & $1$ & & & & & & & \\
$2$ & $1$ & $1$ & & & & & & \\
$3$ & $1$ & $2$ & $2$ & & & & & \\
$4$ & $1$ & $3$ & $5$ & $4$ & & & & \\
$5$ & $1$ & $4$ & $9$ & $12$ & $8$ & & & \\
$6$ & $1$ & $5$ & $14$ & $25$ & $28$ & $16$ & & \\
$7$ & $1$ & $6$ & $20$ & $44$ & $66$ & $64$ & $32$ & \\
$8$ & $1$ & $7$ & $27$ & $70$ & $129$ & $168$ & $144$ & $64$\\
\hline
\end{tabular}
\caption[The number of boolean permutations of each length in $\mathfrak{S}_1, \ldots, \mathfrak{S}_8$.]{The number of boolean permutations of each length in $\mathfrak{S}_1, \ldots, \mathfrak{S}_8$.  Missing table entries are equal to $0$.}\label{L(n,k) table}
\end{table}

It is interesting to note that $B(w)$ is boolean if and only if it is a lattice.  The ideal $B(w)$ is a lattice if and only if all of the $R$-polynomials are of the form $(q-1)^{\ell(y) - \ell(x)}$, as discussed by Brenti in \cite{brenti}.  Moreover, Brenti shows that this is equivalent to all of the Kazhdan-Lusztig polynomials equaling the $g$-polynomials of the duals of the corresponding subintervals.  The $g$-polynomials are defined in \cite{stanleyhvectors}, and their coefficients are the toric $g$-vectors.

Subsequent to Theorem~\ref{boolean thm}, it is natural to ask the following questions.  What can be said about the principal order ideal of permutations with exactly one occurrence of exactly one of the patterns $321$ or $3412$?  (Note that these have reduced decompositions in which exactly one letter is repeated, and it appears exactly twice.)  In particular, how big are these ideals?  These questions are answered below.  Generalizations allowing more occurrences of $321$ and $3412$ are not treated here.

\begin{figure}[htbp]
\centering
\epsfig{file=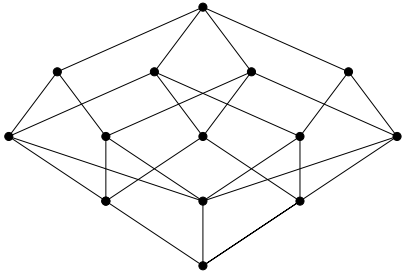, width=2.08in}
\caption{The principal order ideal $B(3412)$.}\label{B(3412) fig}
\end{figure}

\begin{itemize}
\item Suppose that $w$ has exactly one $321$-pattern and is $3412$-avoiding.  Let $\ell = \ell(w)$.  There exists $i_1 \cdots i_{\ell} \in R(w)$ such that $i_j=i_{j+2}$ for some $j$, and there are no repeated letters besides $i_j$ and $i_{j+2}$.  The subword property dictates the poset $B(w)$ as follows.  Consider the poset $B_{\ell}$ of subsets of $[\ell]$ ordered by set inclusion.  Delete all elements of the poset containing $\{j,j+2\}$ but not $j+1$, and identify all elements of the poset containing $j$ but not $\{j+1,j+2\}$ with those that interchange the roles of $j+2$ and $j$.  The resulting poset is isomorphic to $B(w)$, and $|B(w)| = 3 \cdot 2^{\ell - 2}$.\\
Figure~\ref{B(321) fig} depicts $B(w)$ for the simplest such permutation, $w = 321$.
\begin{figure}[htbp]
\centering
\epsfig{file=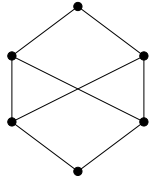, width=.74in}
\caption{The principal order ideal $B(321)$.}\label{B(321) fig}
\end{figure}
\item Suppose that $w$ has exactly one $3412$-pattern and is $321$-avoiding.  Let $\ell = \ell(w)$.  There exists $i_1 \cdots i_{\ell} \in R(w)$ such that $i_j i_{j+1} i_{j+2} i_{j+3} = (2132)^M$ for some $M$, and there are no repeated letters besides $i_j$ and $i_{j+3}$.  Again, consider the poset $B_{\ell}$ of subsets of $[\ell]$ ordered by set inclusion.  Delete all elements of the poset containing $\{j,j+3\}$ but not $\{j+1,j+2\}$, and identify those elements of the poset that contain $j$ but not $\{j+1,j+2,j+3\}$ with those that interchange the roles of $j+3$ and $j$.  The resulting poset is isomorphic to $B(w)$, and $|B(w)| = 7 \cdot 2^{\ell - 3}$.\\
Figure~\ref{B(3412) fig} depicts $B(w)$ for the simplest such permutation, $w = 3412$.
\end{itemize}

\section{Principal order ideals isomorphic to a power of $B(w_0^{(k)})$}\label{power section}

The previous section characterized all permutations for which $B(w)$ is boolean, where a boolean poset is one which is isomorphic to some power of $B(21)$.  This section generalizes the previous work by describing the permutations for which $B(w)$ is isomorphic to a power of $B( w_0^{(k)})$ for $k \ge 3$.

\begin{defn}
Let $k \ge 3$ be an integer and $w \in \mathfrak{S}_n$ be a permutation.  If $B(w) \cong B(w_0^{(k)})^r$ for some $r$, then $w$ is a \emph{power permutation}.
\end{defn}

As in the previous section, the characterization of power permutations is in terms of patterns, although not in quite the same way as Theorem~\ref{boolean thm}.  A few preliminaries are necessary before this characterization can be stated.

\begin{prop}\label{B(k...1) prop}
For $x,y \in \mathfrak{S}_n$, suppose that $[x,y] \cong B(w_0^{(k)})$ for some $k$.  Then
there exist $\bm{i} \in R(x)$ and $\bm{j} \in R(y)$ such that $\bm{i}$ is obtained by deleting a factor from $\bm{j}$ which is the shift of an element of $R(w_0^{(k)})$.
\end{prop}

\begin{proof}
Let $\bm{i} \in R(x)$ and $\bm{j} \in R(y)$ be, by the subword property, such that $\bm{i}$ is a subword of $\bm{j}$.  Consider the multiset $S$ of $\binom{k}{2}$ letters deleted from $\bm{j}$ to form $\bm{i}$.  Because $[x,y] \cong B(w_0^{(k)})$, this $S$ contains $k-1$ distinct letters.

The number of distinct letters in $S$ equals the number of elements covering $x$ in $[x,y]$.  Therefore it must be possible to find $\bm{i}$ and $\bm{j}$ as above so that the factors in $\bm{j}$ formed by $S$ have the property that equal elements of $S$ lie in the same factor.

Given $T$ distinct and consecutive letters, the longest reduced decomposition that can be formed by them has length $\binom{T+1}{2}$.  Observe that
\begin{equation*}
\binom{T_1+1}{2} + \binom{T_2+1}{2} < \binom{T_1+T_2+1}{2}
\end{equation*}
\noindent for $T_1, T_2 > 0$.  Thus, all of $S$ comprises a single factor in $\bm{j}$.
\end{proof}

\begin{prop}\label{k>3 distinct}
If $k \ge 4$, and $\bm{j}^{M_1} \cdots \bm{j}^{M_r} \in R(w)$ for $\bm{j} \in R(w_0^{(k)})$, then the $M_i$s are distinct.
\end{prop}

\begin{proof}
Fix $k \ge 4$.  It is straightforward to show that if $\bm{j}^{M_1} \cdots \bm{j}^{M_r} \in R(w)$, then 
\begin{equation}\label{still decreasing}
w(x + M_i) = k+1-x + M_i \text{ for } x \in [2,k-1],
\end{equation}
for all $i \in [r]$.  The result follows from equation~\eqref{still decreasing} and the fact that elements of $R(w)$ are reduced.
\end{proof}

The result does not hold for $k=3$ because equation~\eqref{still decreasing} says only that $2+M_i$ is fixed by $w$ for all $i$.  For example, $(121)^2(121)^0(121)^4(121)^2 \in R(5274163)$.

\begin{prop}\label{shifts for power}
Fix $k \ge 3$.  Every inversion in $w$ is in exactly one $w_0^{(k)}$-pattern if and only if there exists
\begin{equation}\label{r distinct shifts}
\bm{i} = \bm{j}^{M_1} \cdots \bm{j}^{M_r} \in R(w)
\end{equation}
\noindent for $\bm{j} \in R(w_0^{(k)})$, where the $M_i$s are distinct.  (Consequently there are $r$ occurrences of the pattern $w_0^{(k)}$ in $w$.)  Because $\bm{i}$ is reduced, $|M_i - M_j| \ge k-1$ for all $i \neq j$.
\end{prop}

\begin{proof}
Fix $k \ge 3$.  The result is straightforward for permutations with zero or one $w_0^{(k)}$-pattern.  Suppose that $w \in \mathfrak{S}_n$ has $r>1$ occurrences of $w_0^{(k)}$, and that every inversion in $w$ is in exactly one $w_0^{(k)}$-pattern.  At least one of these patterns occurs in consecutive positions.  Therefore, for some $M$, there exists
\begin{equation*}
w' \eqdef w \cdot \big(1\cdots M (k+M)(k-1+M) \cdots (1+M)(k+1+M) \cdots n\big)
\end{equation*}
\noindent where every inversion in $w'$ is in exactly one $w_0^{(k)}$-pattern, and there are $r-1$ such patterns.  Thus, by induction, there exists $\bm{j}^{M_1} \cdots \bm{j}^{M_r} \in R(w)$ for $\bm{j} \in R(w_0^{(k)})$.  If $k>3$, then this direction of the proof is complete by Proposition~\ref{k>3 distinct}.  If $k=3$ and the $M_i$s are not distinct, then the permutation $w$ would necessarily have a $4312$-, $4231$-, or $3421$-pattern, which contradicts the original hypothesis.

Conversely, suppose that $w$ has a reduced decomposition as in equation~\eqref{r distinct shifts}, where the $M_i$s are distinct.  Consider applying the braid relations to $\bm{i}$.  It is impossible to get a factor equal to the shift of an element of $R(w_0^{(k+1)})$.  Likewise any shift of an element of $R(4312)$, $R(4231)$, or $R(3421)$ can be extended to a factor that is the shift of an element of $R(4321)$.  Therefore the vexillary characterization implies that $w$ avoids $w_0^{(k+1)}$, $4312$, $4231$, and $3421$.  Consequently, every inversion in $w$ is in exactly one $w_0^{(k)}$-pattern.
\end{proof}

With this groundwork, the main theorem of the section can now be stated.

\begin{thm}\label{power thm}
The permutation $w$ is a power permutation if and only if every inversion in $w$ is in exactly one $w_0^{(k)}$-pattern for some fixed $k \ge 3$.
\end{thm}

\begin{proof}
Fix $k \ge 3$ and  $w \in \mathfrak{S}_n$.  First suppose that every inversion in $w$ is in exactly one $w_0^{(k)}$-pattern, and that $w$ contains $R$ distinct occurrences of the pattern $w_0^{(k)}$.  ``Undoing'' inversions in one of these patterns does not alter the other patterns.  Consequently $B(w) \cong B(w_0^{(k)})^R$.

For the other direction of the proof, suppose that $B(w) \cong B(w_0^{(k)})^R$, and proceed by induction on $R$.  The case $R = 0$ is trivial, and the case $R=1$ was  considered in Proposition~\ref{B(k...1) prop}.  Suppose that the theorem holds for permutations whose principal order ideals are isomorphic to $B(w_0^{(k)})^r$, for all $r \in [0,R)$.

There are distinct permutations $w_1, \ldots, w_R$, each less than $w$, with 
\begin{equation*}
B(w_1) \cong \cdots \cong B(w_R) \cong B(w_0^{(k)})^{R-1}.
\end{equation*}
\noindent By Proposition~\ref{shifts for power} and the inductive hypothesis, each of these $R$ permutations has a reduced decomposition $\bm{j}^{M_1^h} \cdots \bm{j}^{M_{R-1}^h} \in R(w_h)$, where the $M_i^h$s are distinct.

The interval $[w_h,w]$ is isomorphic to $B(w_0^{(k)})$ for all $h$, so Proposition~\ref{B(k...1) prop} indicates that $w$ has a reduced decomposition $\bm{j}^{M_1} \cdots \bm{j}^{M_R}$.  The distinct permutations $w_1, \ldots, w_R$ each satisfy the induction hypothesis.  Hence the $M_i$s are distinct, and Proposition~\ref{shifts for power} completes the proof.
\end{proof}

The following corollary is stated in the language of \cite{rdpp}, where $X(w)$ is Elnitsky's polygon (defined in \cite{elnitsky}) and all the polygons have unit sides.

\begin{cor}
If $w$ is a power permutation, then there is a zonotopal tiling of $X(w)$ consisting entirely of $2k$-gons for some $k \ge 3$.  The converse is true if $k \ge 4$.
\end{cor}

Theorem~\ref{power thm} gives a concise description of power permutations, again in terms of patterns.  Although the flavor of this description differs from that of Theorem~\ref{boolean thm}, the prominent role of patterns in the power permutation characterization is immediately apparent.  It is clear that Theorem~\ref{power thm} must be restricted to $k \ge 3$, while the $k=2$ case is treated in Theorem~\ref{boolean thm}, because to say that ``every inversion in $w$ is in exactly one $21$-pattern'' provides no information.

It is instructive to consider what it means for every inversion in $w \in \mathfrak{S}_n$ to be in exactly one $w_0^{(k)}$-pattern.  The following facts are straightforward to show.
\begin{itemize}
\item Distinct occurrences of $w_0^{(k)}$ are disjoint or share exactly one entry.
\item If two occurrences of $w_0^{(k)}$ intersect, then either $\langle k \rangle_1 = \langle k \rangle_2$ or $\langle 1 \rangle_1 = \langle 1 \rangle_2$.
\item Without loss of generality, all values in $\langle w_0^{(k)} \rangle_1$ are at least as large as all values in $\langle w_0^{(k)} \rangle_2$.  The non-shared values in $\langle w_0^{(k)} \rangle_1$ all occur to the right of the non-shared values in $\langle w_0^{(k)} \rangle_2$.
\item If $m$ is not in any $w_0^{(k)}$-pattern, then $w(m) = m$ and $m$ is not in any element of $R(w)$.  Also, $w(1) \cdots w(m-1) \in \mathfrak{S}_{m-1}$ and $\big(w(m+1) - m\big)\cdots \big(w(n) - m\big) \in \mathfrak{S}_{n-m}$ are power permutations with the same parameter $k$.
\item The values $\langle k-1 \rangle, \ldots, \langle 2\rangle$ occur consecutively in $w$.
\end{itemize}

\begin{example}
$521436 \in \mathfrak{S}_6$ and $432159876 \in \mathfrak{S}_9$ are both power permutations.
\end{example}

\section{Patterns and order ideals}\label{pattern ideal section}

The previous sections considered principal order ideals in the Bruhat order of the symmetric group.  This section examines order ideals that are not necessarily principal.  The following questions are completely answered.

\begin{enumerate}
\item For what $p \in \mathfrak{S}_k$, where $k \ge 3$, is the set
\begin{equation*}
S_n\{p\} = \{w \in \mathfrak{S}_n: w \text{ is } p\text{-avoiding}\}
\end{equation*}
\noindent a nonempty order ideal, for some $n > k$?
\item For what $p \in \mathfrak{S}_k$ and $q \in \mathfrak{S}_l$, where $k,l \ge 3$, is the set
\begin{equation*}
S_n\{p,q\} = \{w \in \mathfrak{S}_n: w \text{ is } p \text{- and } q\text{-avoiding}\}
\end{equation*}
\noindent a nonempty order ideal for some $n \ge k,l$?
\end{enumerate}

\noindent The restrictions on $n$, $k$, and $l$ eliminate trivial cases.

Somewhat surprisingly, very few patterns that answer the above questions.

\begin{thm}\label{S_n{p}}
For $k \ge 3$, there is no permutation $p \in \mathfrak{S}_k$ for which there exists $n > k$ such that the set $S_n\{p\}$ is an order ideal.
\end{thm}

\begin{proof}
If $w_0^{(n)}$ avoids $p$, then $S_n\{p\}$ is not an order ideal: $p(1) \cdots p(k) (k+1) \cdots n$ is less than $w_0^{(n)}$ and not in $S_n\{p\}$.  Thus $p = w_0^{(k)}$.

Let $w = k(k+1) (k-1) (k-2) \cdots 4312 (k+2)(k+3) \cdots n \in S_n\{p\}$.
\begin{eqnarray*}
w &\gtrdot& k(k+1) (k-1) (k-2) \cdots 4132 (k+2)(k+3) \cdots n\\
&\gtrdot&  k(k+1) (k-1) (k-2) \cdots 1432 (k+2)(k+3) \cdots n\\
&\gtrdot& \cdots \\
&\gtrdot&  k(k+1) 1 (k-1) (k-2) \cdots 432 (k+2)(k+3) \cdots n\\
&\gtrdot& 1(k+1) k (k-1) (k-2) \cdots 432 (k+2)(k+3) \cdots n =: v.
\end{eqnarray*}
\noindent Because $v$ has a $p$-pattern, the set $S_n\{p\}$ is not an order ideal for any $n > k$.
\end{proof}

The set $S_n\{321,3412\}$ of boolean permutations is an order ideal.  Thus there are permutations $p$ and $q$ for which the set $S_n\{p,q\}$ is an order ideal.

\begin{thm}\label{S_n{p,q}}
Let $p \in \mathfrak{S}_k$ and $q \in \mathfrak{S}_l$ for $k,l \ge 3$.  The only times when $S_n\{p,q\}$ is a nonempty order ideal for some $n \ge k,l$ are $S_n\{321,3412\}$, $S_n\{321,231\}$, and $S_n\{321,312\}$.  These sets are order ideals for all $n \ge 4$
\end{thm}

\begin{proof}
As in the previous proof, it can be assumed that $p = w_0^{(k)}$.

Suppose that $S_n\{p,q\}$ is a nonempty order ideal for some $n \ge k,l$.  Then the following permutations cannot be in $S_n\{p,q\}$, because they are all larger in the Bruhat order than a permutation containing a $p$-pattern.
\begin{eqnarray*}
&k \cdots 3(k+1)12 (k+2)\cdots n&\\
&k(k+1)1(k-1) \cdots 32 (k+2) \cdots n&\\
&1 \cdots (n-k-1) (n-1) \cdots (n-k+2)n(n-k)(n-k+1)&\\
&1 \cdots (n-k-1) (n-1)n(n-k)(n-2) \cdots (n-k+2)(n-k+1)&
\end{eqnarray*}

These all avoid $p$, so they must contain $q$.  The only patterns in all of these permutations are $\{312, 231, 3412, 12 \cdots (l-1)l\}$.  If $q = 12\cdots (l-1)l$ and $S_n\{p,q\}$ is nonempty, then it is not an order ideal because every element in $S_n\{p,q\}$ is greater than $12\cdots n \notin S_n\{p,q\}$.  Therefore $q \in \{312,231,3412\}$.

Suppose that $k>3$.  If $q \in \{231,312\}$, then $u = 32145\cdots n \in S_n\{p,q\}$.   However, $u > q(1)q(2)q(3)45\cdots n \notin S_n\{p,q\}$.  Similarly, $v = 342156\cdots n \in S_n\{p,3412\}$, but $v > 341256\cdots n \notin S_n\{p,3412\}$.  Thus $k=3$ if $S_n\{p,q\}$ is to be an order ideal.

By the vexillary characterization and Theorem~\ref{boolean thm}, the set $S_n\{321,231\}$ consists of those permutations that have reduced decompositions $i_1 \cdots i_{\ell}$ for $i_1 > \cdots > i_{\ell}$.  Thus $S_n\{321,231\}$ is an order ideal by the subword property.  Similarly, the set $S_n\{321,312\}$ consists of those permutations that have reduced decompositions $i_1 \cdots i_{\ell}$ for $i_1 < \cdots < i_{\ell}$.  Once again, this is an order ideal.  As stated earlier, the set $S_n\{321,3412\}$ of boolean permutations is also an order ideal.
\end{proof}

The elements of $S_n\{321,3412\}$ were enumerated by length in Corollary~\ref{boolean by length cor}, and their rank generating function is equation~\eqref{boolean gen fxn}.  The enumerations for the sets $S_n\{321,231\}$ and $S_n\{321,312\}$ are straightforward.

\begin{cor}\label{enumeration of 231/312}
The number of elements of length $k$ in each of $S_n\{321,231\}$ and $S_n\{321,312\}$ is $\binom{n-1}{k}$.  Consequently, each has rank generating function
\begin{equation*}
\sum_{n,k}\binom{n-1}{k} t^kz^n = \frac{z}{1-(1+t)z}.
\end{equation*}
\end{cor}

\begin{proof}
A length $k$ element in $S_n\{321,231\}$ has a reduced decomposition $i_1\cdots i_k$ where $i_1>\cdots>i_k$.  Therefore, it is uniquely determined by choosing $k$ of the $n-1$ possible letters.  The enumeration for $S_n\{321,312\}$ is analogous.
\end{proof}

In each instance where $S_n\{p,q\}$ is an order ideal, the rank generating function of this subposet is a rational function.  For $S_n\{321,231\}$ and $S_n\{321,312\}$, these order ideals are actually \emph{principal}: $S_n\{321,231\} = B(n12\cdots(n-1))$, and $S_n\{321,312\} = B(23\cdots n1)$.  Results of Lakshmibai and Sandhya (see \cite{lakshmibai}) and Carrell and Peterson (see \cite{carrell}) show that $B(w)$ is rank symmetric if and only if $w$ is $3412$- and $4231$-avoiding, which shows (although it is already clear from Corollary~\ref{enumeration of 231/312}) that $S_n\{321,231\}$ and $S_n\{321,312\}$ are both rank symmetric.

The poset of boolean permutations, $S_n\{321,3412\}$, is simplicial, and its $f$-vector was computed in equation~\eqref{L(n,k) equation}.  In \cite{stanleyfvectors}, Stanley showed that for a given vector $\bm{h}$, there exists a Cohen-Macaulay simplicial poset with $h$-vector equal to $\bm{h}$ if and only if $h_0 = 1$ and $h_i \ge 0$ for all $i$.  The last coordinate of the $h$-vector of $S_n\{321,3412\}$ is $L(n,n-1) - L(n,n-2)$, which is negative for $n>3$ (the only $n$ for which $S_n\{321,3412\}$ is defined).  Thus $S_n\{321,3412\}$ is never Cohen-Macaulay.

\section{Boolean order ideals in the Bruhat order for types $B$ and
$D$}\label{bd bruhat section}

As Section~\ref{boolean section} studied boolean principal order ideals in $\mathfrak{S}_n$, this section does likewise for signed permutations.  Recall that the finite Coxeter groups of types $B$ and $D$ consist of signed permutations, where $\D \subset \B$ is the subset of elements that have an even number of negative signs when written in one-line notation.

\begin{example}
$\mathfrak{S}_2^B = \{12, 21, \ul{1}2, \ul{2}1, 1\ul{2}, 2\ul{1}, \ul{12}, \ul{21}\}$ and $\mathfrak{S}_2^D = \{12, 21, \ul{12}, \ul{21}\}$.
\end{example}

The central object here is the principal order ideal of a signed permutation.

\begin{defn}
Let $W$ be a finite Coxeter group of type $A$, $B$, or $D$.  The element $w \in W$ is \emph{boolean} if $B(w)$ is a boolean poset.
\end{defn}

The following proposition holds for $\B$ and $\D$ as well the symmetric group, and its proof is omitted.

\begin{prop}\label{boolean distinct}
Let $W$ be a finite Coxeter group of type $A$, $B$, or $D$.  An element $w \in W$ is boolean if and only if a reduced decomposition of $w$ has no repeated letters.
\end{prop}

Proposition~\ref{boolean distinct} resembles a result of Fan in \cite{fan} for an arbitrary Weyl group $W$.  Fan showed that if the reduced decompositions of $w \in W$ avoid factors of the form $sts$, then the corresponding Schubert variety $X_w$ is smooth if and only if some (every) reduced decomposition of $w$ contains no repeated letter.

The classifications of the boolean elements in $\B$ and $\D$ rely on Proposition~\ref{boolean distinct}.  In each case, the boolean elements are described and enumerated by length.  As with $\mathfrak{S}_n$, these characterizations are in terms of patterns, although the type $B$ case is more complicated than type $A$, and type $D$ is more complicated still.

\begin{thm}\label{type b boolean}
The signed permutation $w \in \B$ is boolean if and only if $w$ avoids all of the following patterns.
\begin{center}
$\begin{array}{ll}
\ul{12} & \ul{21} \\
321 & 3412\\
32\ul{1} & 34\ul{1}2\\
\ul{3}21 & \ul{3}412\\
1\ul{2} & 3\ul{2}1
\end{array}$
\end{center}
\end{thm}

\begin{proof}
By Proposition~\ref{boolean distinct}, a reduced decomposition of a boolean element contains at most one $0$.  Therefore boolean elements in $\B$ have at most one negative value.  Thus the patterns $\ul{12}$ and $\ul{21}$ must be avoided.  Similarly, $w \in \B$ is boolean if and only if it has a reduced decomposition with one of the following forms:
\begin{enumerate}
\item An ordered subset of $[n-1]$;
\item $0 \ \{\text{an ordered subset of } [n-1]\}$; or
\item $\{\text{an ordered subset of } [n-1]\} \ 0$.
\end{enumerate}

By Theorem~\ref{boolean thm}, a reduced decomposition of $v \in \mathfrak{S}_n$ is an ordered subset of $[n-1]$ if and only if $v$ is $321$- and $3412$-avoiding.  The product $s_0v$ changes the sign of the value $1$, while $vs_0$ changes the sign of the value in the first position.  Therefore, a boolean permutation in $\B$ also avoids $32\ul{1}$, $34\ul{1}2$, $\ul{3}21$, and $\ul{3}412$.

Finally, a negative value can appear in a boolean permutation in $\B$ only if it is $\ul{1}$ or occurs in the first position.  Thus the permutation also avoids $1\ul{2}$ and $3\ul{2}1$.
\end{proof}

Proposition~\ref{boolean distinct} states that $w \in \B$ is boolean if and only if it has a reduced decomposition whose letters are all distinct.  Given previous results, the enumeration of these elements is straightforward.  Each of $\{0,1, \ldots, n-1\}$ can appear at most once in a reduced decomposition of a boolean element, so it is necessary only to understand when two ordered subsets of $\{0,1, \ldots, n-1\}$ correspond to the same permutation.  There is a bijection between pairs of commuting elements in $\{s_0, s_1, \ldots, s_{n-1}\}$ and pairs of commuting elements in $\{s_1, \ldots, s_{n-1}, s_n\}$.  Therefore, the work of enumerating boolean elements in $\B$ was already done in Section~\ref{boolean section}.

\begin{cor}
The number of boolean signed permutations in $\B$ is $F_{2n+1}$.
\end{cor}

\begin{proof}
The number of boolean signed permutations in $\B$ is equal to the number of boolean unsigned permutations in $\mathfrak{S}_{n+1}$, which is $F_{2n+1}$ by Corollary~\ref{number boolean}.
\end{proof}

The previous result was also obtained by Fan in \cite{fan}.

\begin{cor}\label{boolean b by length}
The number of boolean signed permutations in $\B$ of length $k$ is
\begin{equation*}
\sum_{i=1}^k \binom{n+1-i}{k+1-i} \binom{k-1}{i-1},
\end{equation*}
\noindent where the (empty) sum for $k=0$ is defined to be $1$.
\end{cor}

\begin{proof}
The number of boolean signed permutations in $\B$ of length $k$ is equal to the number of boolean unsigned permutations in $\mathfrak{S}_{n+1}$ of length $k$.
\end{proof}

The boolean elements of $\D$ are defined and enumerated below.  As for types $A$ and $B$, this characterization is in terms of patterns avoidance.

\begin{thm}\label{type d boolean}
The signed permutation $w \in \D$ is boolean if and only if $w$ avoids all of the following patterns
\begin{center}
$\begin{array}{llllll}
\ul{123}\ & \ul{132}\ & \ul{213}\ & \ul{231}\ & \ul{312}\ & \ul{321}\\
321 & 3412\\
32\ul{1} & 3\ul{12} & 34\ul{1}2 & 34\ul{21}\\
\ul{3}21 & \ul{23}1 & \ul{3}412 & \ul{43}12\\
1\ul{2} & 3\ul{2}1\\
\ul{3}2\ul{1} & \ul{3}4\ul{1}2
\end{array}$
\end{center}

\noindent Note that some of these patterns have an odd number of negative values.
\end{thm}

\begin{proof}
By Proposition~\ref{boolean distinct}, a reduced decomposition of a boolean element has at most one $1'$.  Therefore boolean elements in $\D$ have at most two negative values, so $\ul{123}$, $\ul{132}$, $\ul{213}$, $\ul{231}$, $\ul{312}$, and $\ul{321}$ must be avoided.  Similarly, $w \in \D$ is boolean if and only if it has a reduced decomposition with one of the following forms:
\begin{enumerate}
\item An ordered subset of $[n-1]$;
\item $1' \ \{\text{an ordered subset of } [n-1]\}$; or
\item $\{\text{an ordered subset of } [n-1]\} \ 1'$.
\end{enumerate}

By Theorem~\ref{boolean thm}, a reduced decomposition of $v \in \mathfrak{S}_n$ is an ordered subset of $[n-1]$ if and only if $v$ is $321$- and $3412$-avoiding.  The product $s_{1'}v$ maps the value $1$ to $\ul{2}$ and the value $2$ to $\ul{1}$, while $vs_{1'} = \ul{v(2)v(1)}v(3)\cdots v(n)$.  Therefore, a boolean permutation in $\D$ also avoids $32\ul{1}$, $3\ul{12}$, $34\ul{1}2$, $34\ul{21}$, $\ul{3}21$, $\ul{23}1$, $\ul{3}412$, and $\ul{43}12$.

Finally, since negative values in a boolean permutation in $\D$ can only appear either as $\ul{1}$ and $\ul{2}$ or in the first two positions, the permutation must also avoid the patterns $1\ul{2}$, $3\ul{2}1$, $\ul{3}2\ul{1}$, and $\ul{3}4\ul{1}2$.
\end{proof}

As in types $A$ and $B$, the boolean elements in type $D$ can be enumerated, although this enumeration is not as simple to state as in the other types.  Fan computed these values in \cite{fan}, with the following results.

\begin{cor}[Fan]
For $n \ge 4$, the number of boolean elements in $\D$ is
\begin{equation*}
\frac{13-4b}{a^2(a-b)}a^n + \frac{13-4a}{b^2(b-a)}b^n,
\end{equation*}
\noindent where $a = (3 + \sqrt{5})/2$ and $b = (3 - \sqrt{5})/2$.
\end{cor}

\begin{cor}\label{boolean d by length}
For $n>1$, the number of boolean elements in $\D$ of length $k$ is
\begin{equation}\label{boolean d by length eq}
L^D(n,k) \eqdef L(n,k) + 2L(n,k-1) - L(n-2,k-1) - L(n-2,k-2),
\end{equation}
\noindent where $L(n,k)$ is as defined previously, and $L(n,k)$ is $0$ for any $(n,k)$ on which it is undefined.  $L^D(1,0) = 1$ and $L^D(1,1) = 0$.
\end{cor}

\begin{proof}
These enumerative results follow from Theorem~\ref{type d boolean} and 
Corollary~\ref{boolean by length cor}.  The subtracted terms in 
equation~\eqref{boolean d by length eq} resolve the overcounting that 
occurs when the reduced decompositions of a boolean element contain $1'$ but not $2$.  The case $n=1$ must be treated separately because the only element in $\mathfrak{S}_1^D$ is the identity.
\end{proof}

For small $n$ and $k$, the values $L^D(n,k)$ are displayed in Table~\ref{type d lengths table}.

\begin{table}[htbp]
\centering
\begin{tabular}{|c|ccccccccc|}
\hline
\rule[-2mm]{0mm}{6mm}$L^D(n,k)$ & $k=0$ & $1$ & $2$ & $3$ & $4$ & $5$ & $6$ & $7$ & $8$\\
\hline
\rule[0mm]{0mm}{4mm}$n=1$ & $1$ & $0$ & & & & & & &\\
$2$ & $1$ & $2$ & $1$ & & & & & &\\
$3$ & $1$ & $3$ & $5$ & $4$ & & & & &\\
$4$ & $1$ & $4$ & $9$ & $13$ & $8$ & & & &\\
$5$ & $1$ & $5$ & $14$ & $26$ & $30$ & $16$ & & &\\
$6$ & $1$ & $6$ & $20$ & $45$ & $69$ & $68$ & $32$ & &\\
$7$ & $1$ & $7$ & $27$ & $71$ & $133$ & $176$ & $152$ & $64$ &\\
$8$ & $1$ & $8$ & $35$ & $105$ & $230$ & $373$ & $436$ & $336$ & $128$\\
\hline
\end{tabular}
\caption[The number of boolean elements of each length in $\mathfrak{S}_1^D, \ldots, \mathfrak{S}_8^D$.]{The number of boolean elements of each length in $\mathfrak{S}_1^D, \ldots, \mathfrak{S}_8^D$.  Missing table entries are equal to $0$.}\label{type d lengths table}
\end{table}

\bibliography{patt-bru}
\bibliographystyle{plain}

\end{document}